\documentclass{article}
\usepackage{graphics}
\usepackage{amssymb}
\usepackage{amsmath}
\usepackage[english]{babel}
\newtheorem{theorem}{Theorem}

\newtheorem{corollary}[theorem]{Corollary}
\newtheorem{definition}[theorem]{Definition}

\newtheorem{lemma}[theorem]{Lemma}

\newtheorem{remark}[theorem]{Remark}
\begin{document}
\title{A semi-classical trace formula at a totally degenerate critical level.\\
\large{\textit{Contributions of local extremum}}.}
\author{Brice Camus.\\
Mathematisches Institut der LMU M\"{u}nchen.\\
Theresienstra\ss{}e 39, 80333 Munich, Germany.}
\maketitle
\begin{abstract}
We study the semi-classical trace formula at a critical energy
level for an $h$-pseudo-differential operator on $\mathbb{R}^{n}$
whose principal symbol has a totally degenerate critical point for
that energy. This problem is studied for a large time behavior and
under the hypothesis that the principal symbol of the operator has
a local extremum at the critical point.
\end{abstract}
\section{Introduction.}
\label{intro} The semi-classical trace formula for a self-adjoint
$h$-pseudo-differential operator $P_{h}$, or more generally
$h$-admissible (see \cite{[Rob]}), studies the asymptotic
behavior, as $h$ goes to 0, of the spectral function:
\begin{equation}
\gamma (E,h,\varphi )=\sum\limits_{|\lambda _{j}(h)-E|\leq
\varepsilon }\varphi (\frac{\lambda _{j}(h)-E}{h}),  \label{Def
trace}
\end{equation}
where the $\lambda _{j}(h)$ are the eigenvalues of $P_{h}$ and
where we suppose that the spectrum is discrete in $[E-\varepsilon
,E+\varepsilon ]$, some sufficient conditions for this are given
below. If $p_0$ is the principal symbol of $P_{h}$ we recall that
an energy $E$ is regular when $\nabla p_0(x,\xi )\neq 0$ on the
energy surface $\Sigma _{E}=\{(x,\xi )\in
T^{\star}\mathbb{R}^{n}\text{ }/\text{ }p_0(x,\xi )=E\}$ and
critical when it is not regular.

A classical property is the existence of a link between the
asymptotics of (\ref{Def trace}), as $h$ tends to 0, and the
closed trajectories of the Hamiltonian flow of $p_0$ on $\Sigma
_{E}$
\begin{equation*}
\lim_{h\rightarrow 0}\gamma (E,h,\varphi )\rightleftharpoons
\{(t,x,\xi)\in\mathbb{R}\times \Sigma _{E} \text{ / }
\Phi_{t}(x,\xi)=(x,\xi)\},
\end{equation*}
where  $\Phi_{t}=\mathrm{exp}(tH_{p_0})$ and $H_{p_0}=\partial
_{\xi }p_0.\partial_{x} -\partial _{x}p_0.\partial_{\xi}$. A
non-exhaustive list of references concerning this subject is
Gutzwiller \cite{GUT}, Balian and Bloch \cite{BB} for the physic
literature and for a mathematical point of view Brummelhuis and
Uribe \cite{BU}, Paul and Uribe \cite{PU}, and more recently
Combescure, Ralston and Robert \cite{CRR}, Petkov and Popov
\cite{P-P}, Charbonnel and Popov \cite{C-P}.

The case of a non-degenerate critical energy of the principal
symbol $p_0(x,\xi )$, that is such that the critical-set
$\mathbb{\frak{C}}(p_0) =\{(x,\xi )\in T^{\ast
}\mathbb{R}^{n}\text{ }/\text{ }dp_0(x,\xi )=0\}$ is a compact
$C^{\infty }$ manifold with a Hessian $d^{2}p_0$ transversely
non-degenerate along this manifold, has been investigated first by
Brummelhuis, Paul and Uribe in \cite{BPU}. They treated this
question for quite general operators but for some ''small times'',
that is it was assumed that 0 is the only period of the linearized
flow in $\mathrm{supp}(\hat{\varphi})$ when it is small. Later,
Khuat-Duy in \cite{KhD} and \cite{KhD1} has obtained the
contributions of the non-zero periods of the linearized flow with
the assumption that $\rm{supp}(\hat{\varphi})$ is compact, but for
Schr\"{o}dinger operators with symbol $\xi ^{2}+V(x)$ and a
non-degenerate potential $V$. Our contribution to this subject was
to compute the contributions of the non-zero periods of the
linearized flow for some more general operators, always with
$\hat{\varphi}$ of compact support and under some geometrical
assumptions on the flow (see \cite{Cam0} or \cite{Cam}).

Basically, the asymptotics of (\ref{Def trace}) can be expressed
in terms of oscillatory integrals whose phases are related to the
classical dynamics of $p_0$ on $\Sigma_E$. For $(x_0,\xi_0)$ a
critical point of $p_0$, it is well known that the relation:
\begin{equation}
\mathrm{Ker}(d_{x,\xi}\Phi_{t}(x_0,\xi_0)-\mathrm{Id})\neq \{ 0\},
\end{equation}
leads to the study of degenerate oscillatory integrals. What is
new here is that we examine the case of a totally degenerate
energy, that is such that the Hessian matrix at our critical point
is zero. Hence, the linearized flow for such a critical point
satisfies $d_{x,\xi }\Phi _{t}(x_{0},\xi _{0})=\mathrm{Id}$, for
all $t\in\mathbb{R}$ and the oscillatory integrals we have to
consider are totally degenerate.

The results obtained here are global in time, that is we only
assume that $\mathrm{supp}(\hat{\varphi})$ is compact. The core of
the proof lies in establishing suitable normal forms for our phase
functions and in a generalization of the stationary phase formula
for these normal forms.

\section{Hypotheses and main result.}
Let $P_{h}=Op_{h}^{w}(p(x,\xi ,h))$ be a $h$-pseudodifferential
operator in the class of the $h$-admissible operators with symbol
$p(x,\xi ,h) \sim \sum h^{j}p_{j}(x,\xi )$, i.e. there exist
sequences $(p_{j})_{j}\in \Sigma _{0}^{m}(T^{\ast }\mathbb{R}
^{n})$ and $(R_{N}(h))_{N}$ such that:
\begin{equation*}
P_{h}=\sum\limits_{j<N}h^{j}p_{j}^{w}(x,hD_{x})+h^{N}R_{N}(h),\text{
} \forall N\in \mathbb{N},
\end{equation*}
where $R_{N}(h)$ is a bounded family of operators on
$L^{2}(\mathbb{R}^{n})$, for $h\leq h_{0}$, and:
\begin{equation*}\Sigma _{0}^{m}(T^{\ast }\mathbb{R}^{n})=\{a:T^{\ast }\mathbb{R}%
^{n}\rightarrow \mathbb{C},\text{ }\sup|\partial _{x}^{\alpha
}\partial _{\xi }^{\beta }a(x,\xi )|<C_{\alpha ,\beta }m(x,\xi
),\text{ }\forall \alpha ,\beta \in \mathbb{N}^{n}\},
\end{equation*}
where $m$ is a tempered weight on $T^{\ast }\mathbb{R}^{n}$. For a
detailed exposition on $h$-admissible operators we refer to the
book of Robert \cite{[Rob]}. Let us note $p_{0}(x,\xi )$ the
principal symbol of $P_{h}$, $p_{1}(x,\xi )$ the sub-principal
symbol and $\Phi _{t}=\textrm{exp} (tH_{p_0}):T^{\ast
}\mathbb{R}^{n}\rightarrow T^{\ast }\mathbb{R}^{n}$, the
Hamiltonian flow of $H_{p_0}=\partial _{\xi }p_0
.\partial_{x}-\partial _{x}p_0 .\partial_{\xi }$.

If $E_c$ is a critical value of $p_0$, we study the asymptotics of
the spectral distribution
\begin{equation}
\gamma (E_{c},h)=\sum\limits_{\lambda _{j}(h)\in \lbrack
E_{c}-\varepsilon ,E_{c}+\varepsilon ]}\varphi (\frac{\lambda
_{j}(h)-E_{c}}{h}), \label{Objet trace}
\end{equation}
under the hypotheses $(H_{1})$ to $(H_{4})$ given below. \medskip

$(H_{1})$\textit{ The symbol of }$P_{h}$ \textit{is real.\\There
exists $\varepsilon_{0}>0$} \textit{such that
}$p_0^{-1}([E_{c}-\varepsilon_{0},E_{c}-\varepsilon_{0}])$\textit{
is compact.}\medskip\newline%
Then, by Theorem 3.13 of \cite{[Rob]} the spectrum $\sigma
(P_{h})\cap [E_{c}-\varepsilon ,E_{c}+\varepsilon ]$ is discrete
and consists in a sequence $\lambda _{1}(h)\leq \lambda
_{2}(h)\leq ...\leq \lambda _{j}(h)$ of eigenvalues of finite
multiplicities, if $\varepsilon<\varepsilon_{0}$ and $h$ is small
enough. To simplify notations we write $z=(x,\xi)$ for any point
of the phase space.\medskip

$(H_{2})$\textit{ On the energy surface }$\Sigma
_{E_{C}}=p_0^{-1}(\{E_{c}\})$\textit{, }$p_0$\textit{ has a unique
critical point }$z_{0}=(x_{0},\xi _{0})$\textit{ and near }$z_{0}$
:

\begin{equation*}
p_0(z)=E_{c}+\sum\limits_{j=k}^{N}\mathfrak{p}_{j}(z)+\mathcal{O}(||(z-z_{0})||^{N+1}),\text{
} k>2,
\end{equation*}
\textit{where the functions }$\mathfrak{p}_{j}$\textit{\ are homogeneous of degree }$j$%
\textit{\ in }$z-z_{0}.$

$(H_{3})$\textit{ We have }$\hat{\varphi}\in C_{0}^{\infty
}(\mathbb{R})$\textit{.}\bigskip

With $(H_{2})$ the oscillatory-integrals we will have to consider
are totally degenerate. Hence, they cannot be treated by the
classical stationary phase method. To solve this problem we impose
the following condition on the symbol:\medskip

$(H_{4})$\textit{\ The critical point }$z_{0}$\textit{\ is a local
extremum of }$p_0$\textit{.}
\begin{remark}
\rm{$\mathbf{(}H_{4}\mathbf{)}$ implies that the first non-zero
homogeneous component $\mathfrak{p}_{k}$ is even and is positive
or negative definite and also that $z_0$ is isolated on
$\Sigma_{E_c}$.}
\end{remark}
Since we are interested in the contribution to the trace formula
of the fixed point $z_{0}$, to understand the new phenomenon, it
suffices to study:
\begin{equation}
\gamma _{z_{0}}(E_{c},h)=\frac{1}{2\pi }\mathrm{Tr}\int\limits_{\mathbb{R}}e^{i%
\frac{tE_{c}}{h}}\hat{\varphi}(t)\psi ^{w}(x,hD_{x})\mathrm{exp}(-\frac{i}{h}%
tP_{h})\Theta (P_{h})dt.
\end{equation}
Here $\Theta $ is a function of localization near the critical
energy surface $\Sigma_{E_c}$ and $\psi \in C_{0}^{\infty
}(T^{\ast }\mathbb{R}^{n})$ has an appropriate support near
$z_{0}$. Rigorous justifications are given in section 3 for the
introduction of $\Theta (P_{h})$ and in section 4 for $\psi
^{w}(x,hD_{x})$. Then, the new contribution to the trace formula
is given by:
\begin{theorem}
\label{Main1}Under hypotheses $(H_{1})$ to $(H_{4})$ we obtain:
\begin{equation*}
\gamma _{z_{0}}(E_{c},h)\sim
h^{\frac{2n}{k}-n}(\sum\limits_{j=0}^{N}\Lambda _{j,k}(\varphi
)h^{\frac{j}{k}}+\mathcal{O}(h^{\frac{N+1}{k}})), \text{ as
}h\rightarrow 0,
\end{equation*}
where the $\Lambda _{j,k}$ are some distributions and the leading
coefficient is given by:
\begin{equation}
\Lambda _{0,k}(\varphi ) =\frac{1}{k} \left\langle\varphi
(t+p_{1}(z_{0})),t_{z_0}^{\frac{2n-k}{k}}\right\rangle
\frac{1}{(2\pi
)^{n}}\int\limits_{\mathbb{S}^{2n-1}}|\mathfrak{p}_{k}(\theta
)|^{-\frac{2n}{k}}d\theta,
\end{equation}
with $t_{z_0}=\mathrm{max}(t,0)$ if $z_0$ is a minimum and
$t_{z_0}=\mathrm{max}(-t,0)$ for a maximum.
\end{theorem}
\begin{remark}
\rm{One can derive a full asymptotic expansion as shows Lemma
\ref{Theo IO 1ere carte}.}
\end{remark}
\section{Oscillatory representation.}
Let be $\varphi \in \mathcal{S}(\mathbb{R})$ with
$\hat{\varphi}\in C_{0}^{\infty }(\mathbb{R})$. We recall that:
\begin{equation*} \gamma (E_{c},h)
=\sum\limits_{\lambda _{j}(h)\in I_{\varepsilon }}\varphi
(\frac{\lambda _{j}(h)-E_{c}}{h}), \text{ } I_{\varepsilon }
=[E_{c}-\varepsilon ,E_{c}+\varepsilon ],
\end{equation*}
with $p_0^{-1}(I_{\varepsilon_{0} })$ compact in $T^{\ast
}\mathbb{R}^{n}$ the spectrum of $P_{h}$ is discrete in
$I_{\varepsilon }$ for $\varepsilon<\varepsilon_{0}$ and $h$ small
enough. Now, we localize near the critical energy $E_{c}$ with a
cut-off function $\Theta \in C_{0}^{\infty }(]E_{c}-\varepsilon
,E_{c}+\varepsilon \lbrack )$, such that $\Theta =1$ near $E_{c}$
and $0\leq \Theta \leq 1$ on $\mathbb{R}$. The associated
decomposition is:
\begin{equation}\label{def S1 S2}
\gamma (E_{c},h) =\gamma _{1}(E_{c},h)+\gamma _{2}(E_{c},h),
\end{equation}
with:
\begin{equation*}
\gamma _{1}(E_{c},h)=\sum\limits_{\lambda _{j}(h)\in
I_{\varepsilon }}(1-\Theta )(\lambda _{j}(h))\varphi
(\frac{\lambda _{j}(h)-E_{c}}{h}),
\end{equation*}
\begin{equation*}
\gamma _{2}(E_{c},h)=\sum\limits_{\lambda _{j}(h)\in
I_{\varepsilon }}\Theta (\lambda _{j}(h))\varphi (\frac{\lambda
_{j}(h)-E_{c}}{h}).
\end{equation*}
The asymptotic behavior of $\gamma _{1}(E_{c},h)$ is given by:
\begin{lemma}
$\gamma _{1}(E_{c},h)=\mathcal{O}(h^{\infty })$ as $h\rightarrow
0$.\label{S1(h)=Tr}
\end{lemma}
\noindent \textit{Proof.} Since $\varphi \in
\mathcal{S}(\mathbb{R})$, $\forall k\in \mathbb{N}$, $\exists
C_{k}$ such that $|x^{k}\varphi (x)|\leq C_{k}$ on $\mathbb{R}$.
By Theorem 3.13 of \cite{[Rob]} the number of eigenvalues $N(h)$
lying in $I_{\varepsilon }\cap \rm{supp}(1-\Theta )$ is of order
$\mathcal{O}(h^{-n})$, for $h$ small enough. This gives the
estimate:
\begin{equation*}
|\gamma _{1}(E_{c},h)|\leq N(h)C_{k}|\frac{\lambda
_{j}(h)-E_{c}}{h}|^{-k}.
\end{equation*}
But on the support of $(1-\Theta )$ we have $|\lambda
_{j}(h)-E_{c}|>\varepsilon _{0}>0$. This leads to:
\begin{equation*}
|\gamma _{1}(E_{c},h)|\leq N(h)C_{N}\varepsilon _{0}^{-k}h^{k}\leq
c_{N}h^{k-n}.
\end{equation*}
Since the property is true for all $k\in \mathbb{N}$, this ends
the proof. \hfill{$\blacksquare $}\medskip

Consequently, for the study of $\gamma (E_{c},h)$ modulo
$\mathcal{O}(h^{\infty })$, we have only to consider the quantity
$\gamma _{2}(E_{c},h)$. By inversion of the Fourier transform we
have:
\begin{equation*}
\Theta (P_{h})\varphi (\frac{P_{h}-E_{c}}{h})=\frac{1}{2\pi} \int\limits_{%
\mathbb{R}}e^{i\frac{tE_{c}}{h}}\hat{\varphi}(t)\mathrm{exp}(-\frac{i}{h}%
tP_{h})\Theta (P_{h})dt.
\end{equation*}
Since the trace of the left hand-side is exactly $\gamma
_{2}(E_{c},h)$, we obtain:
\begin{equation}
\gamma _{2}(E_{c},h)
=\frac{1}{2\pi}\mathrm{Tr}\int\limits_{\mathbb{R}}e^{i\frac{tE_{c}}{h}}
\hat{\varphi}(t)\mathrm{exp}(-\frac{i}{h}tP_{h})\Theta (P_{h})dt,
\label{Trace S2(h)}
\end{equation}
and with Lemma \ref{S1(h)=Tr} this gives:
\begin{equation*}
\gamma (E_{c},h)=\frac{1}{2\pi }\mathrm{Tr}\int\limits_{\mathbb{R}}e^{i%
\frac{tE_{c}}{h}}\hat{\varphi}(t)\mathrm{exp}(-\frac{i}{h}tP_{h})\Theta
(P_{h})dt+\mathcal{O}(h^{\infty }).
\end{equation*}
\begin{remark}\rm{Another interest of this formulation is that, under the
geometrical condition to have a "clean" flow, $\gamma (E_{c},h)$
can be expressed, at the first order, as the composition of two
Fourier integral-operators.}
\end{remark}

Let be $U_{h}(t)=\mathrm{exp}(-\frac{it}{h}P_{h})$ the evolution
operator. For each integer $L$ we can approximate $U_{h}(t)\Theta
(P_{h})$, modulo $\mathcal{O}(h^{L})$, by a Fourier
integral-operator, or FIO, depending on a parameter $h$. To give a
precise formulation of this approximation we recall briefly the
principal notions on FIO. Let be $N\in \mathbb{N}$, $\varphi \in
C^{\infty }(\mathbb{R}^{n}\times \mathbb{R}^{N})$ and $a\in
C_{0}^{\infty }(\mathbb{R}^{n}\times \mathbb{R}^{N})$. An
oscillatory-integral with phase $\varphi $ and amplitude $a$ is
\begin{equation}
I(ae^{\frac{i}{h}\varphi })=(2\pi h)^{-\frac{N}{2}}\int\limits_{\mathbb{R}%
^{N}}a(x,\theta )e^{\frac{i}{h}\varphi (x,\theta )}d\theta .
\end{equation}
To attain our objectives, it suffices to consider amplitudes $a$
with compact support. We suppose, as usually, that $\varphi
=\varphi (x,\theta )$ is a non-degenerate phase function, i.e.
\begin{equation*}
d(\partial _{\theta _{1}}\varphi )\wedge ...\wedge d(\partial
_{\theta _{N}}\varphi )\neq 0\text{ on }
\mathbb{\frak{C}}(\varphi) =
\{(x,\theta )\in \mathbb{R}^{n}\times \mathbb{R}^{N}\text{ / }%
d_{\theta }\varphi (x,\theta )=0\}.
\end{equation*}
This implies that $\mathbb{\frak{C}}(\varphi )$ is a sub-manifold
of class $C^{\infty }$ of $\mathbb{R}^{n}\times \mathbb{R}^{N}$
and that:
\begin{equation*}
i_{\varphi }:\left\{
\begin{array}{l}
\mathbb{\frak{C}}(\varphi )\rightarrow T^{\ast }(\mathbb{R}^{n}), \\
(x,\theta )\mapsto (x,d_{x}\varphi (x,\theta )),
\end{array}
\right.
\end{equation*}
is an immersion. In this situation one say that $\varphi $
parameterizes the Lagrangian manifold $\Lambda _{\varphi
}=i_{\varphi }(\mathbb{\frak{C}}(\varphi ))$. Conversely, if
$\Lambda \subset T^{\ast }\mathbb{R}^{n}$ is a Lagrangian
sub-manifold we can always find locally some non-degenerate phase
functions parameterizing $\Lambda ,$ see e.g. \cite{DUI1}. Now if
$\varphi _{1}$ and $\varphi _{2},$ $\varphi _{j}\in C^{\infty
}(\mathbb{R}^{n}\times \mathbb{R}^{N_{j}}),$ are two
non-degenerate phase functions parameterizing locally the same
Lagrangian manifold, i.e. $\Lambda _{\varphi _{1}}\cap U=\Lambda
_{\varphi _{2}}\cap U$, for an open $U\subset T^{\ast
}\mathbb{R}^{n}$, then it is well-known that there exists a
constant $c\in \mathbb{R}$ such that for all $a_{1}\in
C_{0}^{\infty }(\mathbb{R}^{n}\times \mathbb{R}^{N_{1}})$ with
$i_{\varphi _{1}}(\mathrm{supp}(a_{1})\cap \mathbb{\frak{C}}
(\varphi_{1}))\subset U$, small enough, there exists
a sequence $a_{2,j}\in C_{0}^{\infty }(\mathbb{R}^{n}\times \mathbb{R}%
^{N_{2}})$, $j\in \mathbb{N}$, such that for all $L\in \mathbb{N}$
:
\begin{equation*}
I(a_{1}e^{\frac{i}{h}\varphi _{1}})=e^{\frac{i}{h}c}\sum%
\limits_{j<L}h^{j}I(a_{2,j}e^{\frac{i}{h}\varphi
_{2}})+h^{L}r_{L}(h),
\end{equation*}
with $r_{L}(h)$ uniformly bounded in $L^{2}(\mathbb{R}^{n})$,
since $a_{1}\in C_{0}^{\infty}$, and $c$ comes from $S_{\varphi
_{1}}-S_{\varphi _{2}}=c$ on $\Lambda _{\varphi _{1}}\cap
U=\Lambda _{\varphi _{2}}\cap U$. This recalls that the
fundamental object associated to an oscillatory-integral is the
Lagrangian manifold parameterized by the phase function. The next
definition follows H\"{o}rmander.
\begin{definition}
Let be $\Lambda \subset T^{\ast }\mathbb{R}^{n}$ a Lagrangian
sub-manifold of class $C^{\infty }$. The class
$I(\mathbb{R}^{n},\Lambda )$ of oscillatory functions associated
to $\Lambda $ is given by:
\begin{equation*}
u_{h}\in I(\mathbb{R}^{n},\Lambda )\Leftrightarrow
u_{h}=\sum\limits_{v}I(a_{v}(h)e^{\frac{i}{h}\varphi _{v}}),
\end{equation*}
where the sum is locally finite and the functions $a_{v}(h)$ are
of the form:
\begin{equation*}
a_{v}(h)=\sum\limits_{-d_{v}\leq j<J_{v}}h^{j}a_{j},\text{
}a_{j}\in C^{\infty }(\mathbb{R}^{n}\times \mathbb{R}^{N_{v}}).
\end{equation*}
With $d_{v},J_{v}$ and $N_{v}$ in $\mathbb{N}$ and where $\varphi
_{v}$ is a non degenerate phase function parameterizing $\Lambda
\cap i_{\varphi
_{v}}(\mathrm{supp}(a_{v}(h))),$ where $\mathrm{supp}(a_{v}(h))=\bigcup%
\limits_{j}\mathrm{supp}(a_{v,j}).$
\end{definition}
Now let be $\Lambda \subset T^{\ast }(\mathbb{R}^{n}\times
\mathbb{R}^{l})$ a $C^{\infty }$ Lagrangian sub-manifold of
$T^{\ast }(\mathbb{R}^{n}\times \mathbb{R}^{l})$.
\begin{definition}
The family of operators $(F_{h}):L^{2}(\mathbb{R}^{n})\rightarrow
L^{2}(\mathbb{R}^{l})$, $0<h\leq1$, is a $h$-FIO, associated to
$\Lambda ,$ if there exists two families $\hat{F}_{h}^{(N)}$ and
$R_{h}^{(N)},$ $N\in \mathbb{N}$, such that:
\newline
i) Each $\hat{F}_{h}^{(N)}$ has an integral kernel in
$I(\mathbb{R}^{n}\times \mathbb{R}^{l},\Lambda ).$\newline%
ii) For each $N$ the operator $R_{h}^{(N)}$ is bounded on $L^{2}$
and there exists $C_{N}>0$ such that
$||R_{h}^{(N)}||_{\mathcal{L}(L^{2}
(\mathbb{R}^{n}),L^{2}(\mathbb{R}^{l}))}\leq C_{N},\text{
uniformly in }h,\text{ }0<h\leq 1$.\newline iii) For each
$N:F_{h}=\hat{F}_{h}^{(N)}+h^{N}R_{h}^{(N)}$.
\end{definition}

After these definitions, we recall the theorem that gives the
approximation of $U_{h}(t)\Theta (P_{h})$. Let $\Lambda$ be the
Lagrangian manifold associated to the flow of $p_0$:
\begin{equation*}
\Lambda =\{(t,\tau ,x,\xi ,y,\eta )\in T^{\ast }\mathbb{R}\times T^{\ast }%
\mathbb{R}^{n}\times T^{\ast }\mathbb{R}^{n}:\tau =p_0(x,\xi
),\text{ }(x,\xi )=\Phi _{t}(y,\eta )\}.
\end{equation*}
\begin{theorem}
The operator $U_{h}(t)\Theta (P_{h})$ is an $h$-FIO associated to
$\Lambda$, there exists $U_{\Theta ,h}^{(N)}(t)$ with integral
kernel in $I(\mathbb{R}^{2n+1},\Lambda )$ and $R_{h}^{(N)}(t)$
bounded, with a $L^{2}$-norm uniformly bounded for $0<h\leq 1$ and
$t$ in a compact subset of $\mathbb{R}$, such that
$U_{h}(t)\Theta(P_{h})=U_{\Theta,h}^{(N)}(t)+h^{N}R_{h}^{(N)}(t)$.
\end{theorem}

We refer to Duistermaat \cite{DUI1} for a proof of this theorem.
\begin{remark}
\rm{By a theorem of Helffer and Robert, see e.g. \cite{[Rob]},
Theorem 3.11 and Remark 3.14, $\Theta (P_{h})$
is an $h$-admissible operator with a symbol of compact support in $%
p_0^{-1}(I_{\varepsilon }).$ This allows us to consider only
oscillatory-integrals with compact support.}
\end{remark}

For our goal the following corollary is crucial.
\begin{corollary}
Let be $\Theta _{1}\in C_{0}^{\infty }(\mathbb{R})$ such that
$\Theta _{1}=1$ on $\rm{supp}(\Theta )$ and $\rm{supp}(\Theta
_{1})\subset I_{\varepsilon }$, then $\forall N\in \mathbb{N}$:
\begin{equation*}
\mathrm{Tr}(\Theta (P_{h})\varphi (\frac{P_{h}-E_{c}}{h}))=\frac{1}{2\pi }%
\mathrm{Tr}\int\limits_{\mathbb{R}}\hat{\varphi}(t)e^{\frac{i}{h}%
tE_{c}}U_{\Theta ,h}^{(N)}(t)\Theta
_{1}(P_{h})dt+\mathcal{O}(h^{N}).
\end{equation*}
\end{corollary}
\noindent \textit{Proof.} By cyclicity of the trace, for all $N\in
\mathbb{N}$ we have:
\begin{equation*}
\mathrm{Tr}(\Theta (P_{h})\varphi (\frac{P_{h}-E_{c}}{h}))
=\frac{1}{2\pi}\mathrm{Tr}\int\limits_{\mathbb{R}}\hat{\varphi}(t)e^{
\frac{i}{h}tE_{c}}(U_{\Theta,h}^{(N)}(t)+h^{N}R_{h}^{(N)}(t))
\Theta_{1}(P_{h})dt.
\end{equation*}
But $\Theta _{1}(P_{h})$, and a fortiori $\Theta
_{1}(P_{h})R_{h}^{(N)}(t)$ by the ideal property, is a class-trace
operator (see e.g. \cite{[Rob]}, section 2.5), with norm:
\begin{equation*}
||\Theta _{1}(P_{h})R_{h}^{(N)}(t)||_{\mathrm{Tr}}
\leq||R_{h}^{(N)}(t)|| \text{ }||\Theta_{1}(P_{h})||_{\mathrm{Tr}}
\leq C_{\hat{\varphi}}||\Theta _{1}(P_{h})||_{\mathrm{Tr}},
\end{equation*}
since it is true for all $t\in \rm{supp}(\hat{\varphi})$ the
corollary is proven. \hfill{$\blacksquare$}\medskip

The geometrical approach associated to the H\"{o}rmander class $I(\mathbb{R}%
^{2n+1},\Lambda )$ gives a great degree of freedom in the choice
of the phase function, this will be exploited below. In fact if
$(x_{0},\xi _{0})\in \Lambda $ and if $\varphi =$ $\varphi
(x,\theta )\in C^{\infty }(\mathbb{R}^{n}\times \mathbb{R}^{N})$
parameterizes $\Lambda $ in a neighborhood $U,$ small enough, of
$(x_{0},\xi _{0})$ then for each $u_{h}\in
I(\mathbb{R}^{n},\Lambda )$ and $\chi \in C_{0}^{\infty }(T^{\ast
}\mathbb{R}^{n})$, $\rm{supp}(\chi )\subset U,$ there exists a
sequence of amplitudes $a_{j}=a_{j}(x,\theta )\in C_{0}^{\infty }(\mathbb{R}%
^{n}\times \mathbb{R}^{N})$ such that for all integer $L$:
\begin{equation}
\chi ^{w}(x,hD_{x})u_{h}=\sum\limits_{-d\leq j<L}h^{j}I(a_{j}e^{\frac{i}{h}%
\varphi })+\mathcal{O}(h^{L}).
\end{equation}
We will use this remark with the following result of H\"{o}rmander
(see \cite {HOR1}, tome 4, proposition 25.3.3). Let $(T,\tau
,x_{0},\xi _{0},y_{0},-\eta _{0})\in \Lambda _{\mathrm{flow}}$,
$\eta _{0}\neq 0$, then near this point there exists, after
perhaps a change of local coordinates in $y$ near $y_{0},$ a
function $S(t,x,\eta )$ such that:
\begin{equation}
\phi (t,x,y,\eta )=S(t,x,\eta )-\left\langle y,\eta \right\rangle,
\end{equation}
parameterizes $\Lambda _{\mathrm{flow}}$. In particular this
implies that:
\begin{equation*}
\{(t,\partial _{t}S(t,x,\eta ),x,\partial _{x}S(t,x,\eta
),\partial _{\eta }S(t,x,\eta ),-\eta )\}\subset \Lambda
_{\mathrm{flow}},
\end{equation*}
and that the function $S$ is a generating function of the flow,
i.e.
\begin{equation}
\Phi _{t}(\partial _{\eta }S(t,x,\eta ),\eta ) =(x,\partial
_{x}S(t,x,\eta )). \label{Gene}
\end{equation}
Moreover, $S$ satisfies the Hamilton-Jacobi equation:
\begin{equation*}
\left\{
\begin{array}{c}
\partial _{t}S(t,x,\eta )+ p_0(x,\partial
_{x}S(t,x,\eta ))=0, \\
S(0,x,\xi)=\left\langle x,\xi \right\rangle .
\end{array}
\right.
\end{equation*}
Now, we apply this result with $(x_{0},\xi _{0})=(y_{0},\eta
_{0})$, our unique fixed point of the flow on the energy surface
$\Sigma _{E_{c}}$. If $\xi _{0}=0$ we can replace the operator
$P_{h}$ by $e^{\frac{i}{h}\left\langle x,\xi _{1}\right\rangle }P_{h}e^{-\frac{i}{h}%
\left\langle x,\xi _{1}\right\rangle }$ with $\xi _{1}\neq 0.$
This will not change the spectrum since this new operator has the
symbol $p(x,\xi -\xi _{1},h)$ and the critical point is now
$(x_{0},\xi _{1})$ with $\xi _{1}\neq 0$. Consequently, the
localized trace $\gamma _{2}(E_{c},h)$, defined by Eq.(\ref{def S1
S2}), can be written for all $N\in \mathbb{N}$ and modulo
$\mathcal{O}(h^{N})$ as:
\begin{equation}
\gamma _{2}(E_{c},h)=\sum\limits_{j<N}(2\pi h)^{-d+j}\int\limits_{\mathbb{%
R\times R}^{2n}}e^{\frac{i}{h}(S(t,x,\xi )-\left\langle x,\xi
\right\rangle +tE_{c})}a_{j}(t,x,\xi )\hat{\varphi}(t)dtdxd\xi .
\label{gamma1 OIF}
\end{equation}
To obtain the right power $-d$ of $h$ occurring in Eq.(\ref{gamma1
OIF}) we apply results of Duistermaat \cite{DUI1} (following here
H\"{o}rmander for the FIO, see \cite {HOR2} tome 4, for example)
concerning the order. An $h$-pseudodifferential operator obtained
by Weyl quantization:
\begin{equation*}
(2\pi h)^{-\frac{N}{2}}\int\limits_{\mathbb{R}^{N}}a(\frac{x+y}{2},\xi )e^{%
\frac{i}{h}\left\langle x-y,\xi \right\rangle }d\xi ,
\end{equation*}
is of order 0 w.r.t. $1/h$. Now, since the order of $U_{h}(t)\Theta (P_{h})$ is $-%
\frac{1}{4}$, we find that
\begin{equation}
\psi ^{w}(x,hD_{x})U_{h}(t)\Theta (P_{h})\sim
\sum\limits_{j<N}(2\pi
h)^{-n+j}\int\limits_{\mathbb{R}^{n}}a_{j}(t,x,y,\eta )e^{\frac{i}{h}%
(S(t,x,\eta )-\left\langle y,\eta \right\rangle )}dy,
\label{operateur d'evolution}
\end{equation}
Multiplying Eq.(\ref{operateur d'evolution}) by $\hat{\varphi}%
(t)e^{\frac{i}{h}tE_{c}}$ and passing to the trace we find
(\ref{gamma1 OIF}) with $d=n$, where we write again
$a_{j}(t,x,\eta )$ for $a_{j}(t,x,x,\eta )$.

To each element $u_{h}$ of $I(\mathbb{R}^{n},\Lambda )$ we can
associate a principal symbol $e^{\frac{i}{h}S}\sigma
_{\mathrm{princ}}(u_{h})$, where $S$ is a function on $\Lambda $
such that $\xi dx=dS$ on $\Lambda .$ If
$u_{h}=I(ae^{\frac{i}{h}\varphi })$, then we have $S=S_{\varphi
}=\varphi \circ i_{\varphi }^{-1}$ and $\sigma
_{\mathrm{princ}}(u_{h})$ is a section of $|\Lambda
|^{\frac{1}{2}} \otimes M(\Lambda )$, where $M(\Lambda )$ is the
Maslov vector-bundle of $\Lambda $ and $|\Lambda |^{\frac{1}{2}}$
the bundle of half-densities on $\Lambda$. If $p_{1}$ is the
sub-principal symbol of $P_{h}$, the half-density of the
propagator $U_{h}(t)$ is easily expressed in the global
coordinates $(t,y,\eta )$ on $\Lambda _{\mathrm{flow}}$ via:
\begin{equation}
\exp (i\int\limits_{0}^{t}p_{1}(\Phi _{s}(y,-\eta ))ds)|dtdyd\eta |^{\frac{1%
}{2}}.\label{demi densite}
\end{equation}
This expression is related to the resolution of the first
transport equation for the propagator. For a proof we refer to
Duistermaat and H\"{o}rmander \cite{D-H}.
\section{Classical dynamics near the critical point.}
A critical point of the phase function of Eq.(\ref{gamma1 OIF})
leads to the equations:
\begin{equation*}
\left\{
\begin{array}{c}
E_{c}=-\partial _{t}S(t,x,\xi ), \\
x=\partial _{\xi }S(t,x,\xi ), \\
\xi =\partial _{x}S(t,x,\xi ),
\end{array}
\right. \Leftrightarrow \left\{
\begin{array}{c}
p_0(x,\xi )=E_{c}, \\
\Phi _{t}(x,\xi )=(x,\xi ),
\end{array}
\right.
\end{equation*}
where the right hand side defines a closed trajectory of the flow
inside $\Sigma _{E_{c}}$. Since we are interested in the
contribution of the critical point, we choose a function $\psi \in
C_{0}^{\infty }(T^{\ast }\mathbb{R}^{n})$, with $\psi =1\text{
near }z_{0}$, hence:
\begin{eqnarray*}
\gamma _{2}(E_{c},h) &=&\frac{1}{2\pi }\mathrm{Tr}\int\limits_{\mathbb{R}}e^{i%
\frac{tE_{c}}{h}}\hat{\varphi}(t)\psi ^{w}(x,hD_{x})\mathrm{exp}(-\frac{i}{h}%
tP_{h})\Theta (P_{h})dt \\
&+&\frac{1}{2\pi }\mathrm{Tr}\int\limits_{\mathbb{R}}e^{i\frac{tE_{c}}{h}}\hat{%
\varphi}(t)(1-\psi
^{w}(x,hD_{x}))\mathrm{exp}(-\frac{i}{h}tP_{h})\Theta (P_{h})dt.
\end{eqnarray*}
The contribution on $\rm{supp}(1-\psi)$ is non-singular and can be
treated by the regular trace formula. If $\rm{supp}(\psi)$ is
small enough only points $(t,x,\xi )=(t,z_{0})$, for $t\in
\rm{supp}(\hat{\varphi})$, will give a contribution, since $z_0$
is isolated on $\Sigma_{E_c}$.

Now, we restrict our attention to the contribution of the critical
point. Until further notice, the derivatives $d$ will be taken
with respect to initial conditions. Since $z_0$ is totally
degenerate, we obtain:
\begin{equation}
d\Phi _{t}(z_{0})=\mathrm{exp}(0)=\mathrm{Id},\text{ }\forall t.
\end{equation}
With $(H_{2})$, the next homogeneous components of the flow are
given by:
\begin{equation}d^{j}\Phi
_{t}(z_{0})=0,\text{ }\forall t,\text{ }\forall j\in \{2,..,k-2\}.
\end{equation}
To obtain the next non-zero terms of the Taylor expansion of the
flow, we will use the following technical result:
\begin{lemma}
\label{TheoFormule de récurence du flot}Let be $z_{0}$ an
equilibrium of the $C^{\infty}$ vector field $X$ and $\Phi _{t}$
the flow of $X$. Then for all $m\in \mathbb{N}^{\ast }$, there
exists a polynomial map $P_{m}$, vector valued and of degree at
most $m$, such that:
\begin{equation}
d^{m}\Phi _{t}(z_{0})(z^{m})=d\Phi
_{t}(z_{0})\int\limits_{0}^{t}d\Phi _{-s}(z_{0})P_{m}(d\Phi
_{s}(z_{0})(z),...,d^{m-1}\Phi _{s}(z_{0})(z^{m-1}))ds.
\label{Formule de récurence du flot}
\end{equation}
\end{lemma}
{\it Proof.} We note $x^l\in (\mathbb{R}^n)^l$ the image of $x$
under the diagonal mapping, with the same convention for any
vector. If $f,g$ are two $C^{\infty }$ applications, by the "Faa
di Bruno formula", we have:
\begin{equation*}
d^{m}(fog)(x_{0})(x^{m})=\sum\limits_{\mathrm{p}\in\mathcal{P}(m)}
C(\mathrm{p}) d^{r}f(g(x_{0}))(d^{n_1}g(x_{0})(x^{n_1})
,...,d^{n_r}g(x_{0})(x^{n_r})),
\end{equation*}
where $\mathcal{P}(m)$ is the set of partitions of the integer $m$
and where it is assumed that the partition p of $m$ is given by
$m=n_1+...+n_r$ and $C(\mathrm{p})$ are some universal integers.
Hence, for our fixed point $z_{0}$ we obtain:
\begin{equation*}
d^{m}(Xo\Phi
_{s})(z_{0})(z^{m})=\sum\limits_{\mathrm{p}\in\mathcal{P}(m)}
C(\mathrm{p})d^{r}X(z_{0})(d^{n_1}\Phi _{s}(z_{0})(z^{n_1})
,...,d^{n_r}\Phi _{s}(z_{0})(z^{n_r})).
\end{equation*}
For $Y=(Y_{1},...,Y_{m})$, we can define:
\begin{equation}
P_{m}(Y)=\sum\limits_{\mathrm{p}\in\mathcal{P}(m)}
C(\mathrm{p})d^{r}X(z_{0})(Y_{n_1},...,Y_{n_r})-dX(z_{0})(Y_{m}),
\end{equation}
this leads to the differential equation, operator valued:
\begin{equation*}
\frac{d}{ds}(d^{m}\Phi _{s}(z_{0}))(z^{m})=dX(z_{0})d^{m}\Phi
_{s}(z_{0})(z^{m})+P_{m}(d\Phi _{s}(z_{0})(z),...,d^{m-1}\Phi
_{s}(z_{0})(z^{m-1})).
\end{equation*}
With the initial condition $d^{m}\Phi _{0}(z_{0})=0,$ the solution
is given by Eq.(\ref{Formule de récurence du flot}).\medskip
\hfill{$\blacksquare$}

Since $d\Phi_{t}(z_{0})=\mathrm{Id}$, $\forall t$, the first
non-zero term of the Taylor expansion is:
\begin{equation}
d^{k-1}\Phi_{t}(z_{0})(z^{k-1})=\int\limits_{0}^{t}
d^{k-1}H_{p_0}(z_{0})(z^{k-1})ds=td^{k-1}H_{\mathfrak{p}_{k}}(z_0)(z^{k-1}),
\label{derive ordre k-1 du flot}
\end{equation}
where the last identity is obtained using that
$d^{k-1}H_{p_0}(z_{0})=d^{k-1}H_{\mathfrak{p}_{k}}(z_0)$.
Moreover, with $d^{2}p_0(z_0)=0$, for the next term Lemma
\ref{TheoFormule de récurence du flot} gives:
\begin{equation*}
d^{k}\Phi_{t}(z_{0})(z^{k})=\int\limits_{0}^{t}d^{k}H_{p_0}(z_{0})(z^{k})ds
=td^{k}H_{\mathfrak{p}_{k+1}}(z_{0})(z^{k}).
\end{equation*}
\begin{remark}\label{struc flow}
\rm{For the derivatives $d^{j}\Phi_{t}(z_0)$, with $j>k$ there is
two different kind of terms, namely:
\begin{equation*}
\int\limits_{0}^{t}d^{j}H_{p_0}(z_{0})(z^{j})ds
=td^{j}H_{p_0}(z_{0})(z^{j})=\mathcal{O}(||z||^{j}),
\end{equation*}
and terms involving powers of $t$. For example, we have:
\begin{gather*}
\int\limits_{0}^{t}d^{j+2-k}H_{p_0}(z_{0})(z^{j+1-k},d^{k-1}\Phi_{t}(z_{0})(z^{k-1}))\\
=\frac{t^{2}}{2}d^{j+2-k}H_{p_0}(z_{0})(z^{j+1-k},d^{k-1}H_{\mathfrak{p}_{k}}(z_0)(z^{k-1})).
\end{gather*}
This term is simultaneously $\mathcal{O}(t^2)$ and
$\mathcal{O}(||z||^j)$ near $(0,z_0)$. A similar result holds for
other terms, which are $\mathcal{O}(t^d)$ for $d\geq2$, by an easy
recurrence.}
\end{remark}
\begin{lemma}\label{structure phase} Near $z_0$, here supposed to be 0 to simplify, we
have:
\begin{equation}
S(t,x,\xi)-\left\langle x,\xi\right\rangle+tE_{c}=
-t(\mathfrak{p}_{k}(x,\xi)+R_{k+1}(x,\xi)+tG_{k+1}(t,x,\xi )),
\label{forme phase}
\end{equation}
where $R_{k+1}(x,\xi ) =\mathcal{O}(||(x,\xi )||^{k+1})$ and
$G_{k+1}(t,x,\xi)=\mathcal{O}(||(x,\xi )||^{k+1})$, uniformly for
$t$ in a compact subset of $\mathbb{R}$.
\end{lemma}
\textit{Proof.} By Taylor and under hypothesis $(H_2)$ we obtain:
\begin{equation}\label{Taylor-flow}
\Phi _{t}(x,\xi ) =(x,\xi )+\frac{1}{(k-1)!}d^{k-1}\Phi
_{t}(0)(z^{k-1})+\mathcal{O}(||z||^{k}).
\end{equation}
Now, we search our local generating function as:
\begin{equation*}
S(t,x,\xi)=-tE_c +\left\langle x,\xi
\right\rangle+\sum\limits_{j=3}^N S_j(t,x,\xi)
+\mathcal{O}(||(x,\xi)||^{N+1}),
\end{equation*}
where the functions $S_j$ are time dependant and homogeneous of
degree $j$ w.r.t. $(x,\xi)$. With the implicit relation:
$\Phi_t(\partial_{\xi}S(t,x,\xi),\xi)=(x,\partial_{x}S(t,x,\xi))$
and using Eq.(\ref{Taylor-flow}) we have:
\begin{equation*}
S(t,x,\xi)=-tE_c+\left\langle x,\xi
\right\rangle+S_{k}(t,x,\xi)+\mathcal{O}(||(x,\xi )||^{k+1}),
\end{equation*}
and comparing terms of degree $k-1$ gives:
\begin{equation*}
J\nabla S_k(t,x,\xi)=-\frac{1}{(k-1)!}d^{k-1}\Phi
_{t}(0)((x,\xi)^{k-1}),
\end{equation*}
where $J$ is the matrix of the usual symplectic form. By
homogeneity and with Eq.(\ref{derive ordre k-1 du flot}) we obtain
that:
\begin{equation*}
S_k(t,x,\xi)=\frac{1}{k!}\left\langle
(x,\xi),tJd^{k-1}H_{\mathfrak{p}_k}(x,\xi)^{k-1}
\right\rangle=-t\mathfrak{p}_{k}(x,\xi).
\end{equation*}
It remains now to treat the remainder. Since $S(0,x,\xi)=
\left\langle x,\xi \right\rangle $, we have:
\begin{equation*}
S(t,x,\xi)-\left\langle x,\xi \right\rangle=tF(t,x,\xi),
\end{equation*}
with $F$ smooth in a neighborhood of $(x,\xi)=0$. Now, the
Hamilton-Jacobi equation imposes that $F(0,x,\xi)=-p_0(x,\xi)$ and
we obtain:
\begin{equation*}
R_{k+1}(x,\xi)=p_0(x,\xi)-E_c-\mathfrak{p}_k(x,\xi)=\mathcal{O}(||(x,\xi)||^{k+1}).
\end{equation*}
Finally, the time dependant remainder can be written:
\begin{equation*}
S(t,x,\xi)-S(0,x,\xi)-t\partial_t S(0,x,\xi)=\mathcal{O}(t^2),
\end{equation*}
since by construction this term is of order
$\mathcal{O}(||(x,\xi)||^{k+1})$ we get the result.
$\hfill{\blacksquare}$
\section{Normal forms of the phase function.}
Since the contribution we study is local, we can work with some
coordinates and identify locally $T^{\ast }\mathbb{R}^{n}$ with
$\mathbb{R}^{2n}$ near the critical point. We define:
\begin{equation}
\Psi(t,z)=\Psi(t,x,\xi)=S(t,x,\xi )-\left\langle x,\xi
\right\rangle+tE_{c},\text{ }z=(x,\xi)\in \mathbb{R}^{2n}.
\label{defphase}
\end{equation}
\begin{lemma}\label{FN1}
If $P_h$ satisfies conditions $(H_{2})$ and $(H_{4})$ then, in a
neighborhood of $z=z_0$, there exists local coordinates $\chi$
such that $ \Psi(t,z)\simeq \chi _{0}\chi _{1}^{k}$ if $z_0$ is a
maximum and $ \Psi(t,z)\simeq -\chi _{0}\chi _{1}^{k}$ if $z_0$ is
minimum.
\end{lemma}
\noindent\textit{Proof.} We can here assume that $z_{0}$ is the
origin and we use polar coordinates $z=(r,\theta ),$ $\theta \in
\mathbb{S}^{2n-1}(\mathbb{R})$. With Lemma \ref{structure phase},
near the critical point we have:
\begin{equation*}
\Psi(t,z)\simeq \Psi(t,r\theta)=-tr^{k}(\mathfrak{p}_{k}(\theta
)+rR_{k+1}(\theta )+tG_{k+1}(t,r\theta )),
\end{equation*}
where $\mathfrak{p}_{k}(\theta )$ is the restriction of
$\mathfrak{p}_{k}$ on $\mathbb{S}^{2n-1}$. We define new
coordinates:
\begin{gather*}
(\chi _{0},\chi _{2},...,\chi _{2n})(t,r,\theta) =(t,\theta _{1},...,\theta _{2n-1}), \\
\chi _{1}(t,r,\theta) = r|\mathfrak{p}_{k}(\theta
)+rR_{k+1}(\theta )+tG_{k+1}(t,r\theta )|^{\frac{1}{k}}.
\end{gather*}
In these coordinates the phase becomes $-\chi _{0}\chi _{1}^{k}$
for a minimum and $\chi _{0}\chi _{1}^{k}$ for a maximum. Now, we
observe that:
\begin{equation*}
\frac{\partial\chi_1}{\partial
r}(t,0,\theta)=|\mathfrak{p}_k(\theta)|^{\frac{1}{k}},\text { }
\forall t,
\end{equation*}
hence, the corresponding Jacobian satisfies the relation:
\begin{equation*}
|J\chi |(t,0,\theta )=|\mathfrak{p}_{k}(\theta
)|^{\frac{1}{k}}\neq 0,\text{ } \forall t,
\end{equation*}
and this defines a local system of coordinates near $z_0$ for all
$t$.\hfill{$\blacksquare$}

This change of coordinates leads to:
\begin{equation*}
\int\limits_{\mathbb{R\times R}_{+}\times
\mathbb{S}^{2n-1}}e^{\frac{i}{h}\Psi (t,r,\theta )}a(t,r\theta
)r^{2n-1}dtdrd\theta =\int e^{\pm \frac{i}{h}\chi _{0}\chi
_{1}^{k}}A(\chi _{0},\chi _{1})d\chi _{0}d\chi _{1},
\end{equation*}
with a new amplitude $A$ defined by:
\begin{equation}
A(\chi _{0},\chi _{1})=\int \chi ^{\ast }(a(t,r\theta
)r^{2n-1}|J\chi |)d\chi _{2}...d\chi _{2n}.
\end{equation}
\begin{remark}\label{shift}
\rm{Since $\chi _{1}(t,r,\theta ) =r|\mathfrak{p}_{k}(\theta
)+rR_{k+1}(\theta )+tG_{k+1}(t,r\theta )|^{\frac{1}{k}}$, our new
amplitude satisfies $A(\chi _{0},\chi _{1})=\mathcal{O}(\chi
_{1}^{2n-1})$, near $\chi _{1}=0$. This will play a major role
since distributions of Lemma \ref{Theo IO 1ere carte} below are
supported in $\{ \chi_1=0\}$.}
\end{remark}

We end this section with a lemma on asymptotics of oscillatory
integrals.
\begin{lemma}
\label{Theo IO 1ere carte}There exists a sequence $(c_{j})_{j}$ of
distributions, whose support is contained in the set $\{\chi
_{1}=0\},$ such that for all function $a\in C_{0}^{\infty }(\mathbb{R}%
_{+}$ $\times \mathbb{R})$:
\begin{equation}
\int\limits_{0}^{\infty }(\int\limits_{\mathbb{R}}e^{i\lambda \chi
_{0}\chi _{1}^{k}}a(\chi _{0},\chi _{1})d\chi _{0})d\chi _{1}\sim
\sum\limits_{j=0}^{\infty }\lambda ^{-\frac{j+1}{k}}c_{j}(a),
\end{equation}
asymptotically for $\lambda \rightarrow \infty ,$ with
\begin{equation*}
c_{j}=\frac{1}{k}\frac{1}{j!}(\mathcal{F}(x_{-}^{\frac{j+1-k}{k}})(\chi_0)
\otimes\delta _{0}^{(j)}(\chi_1)),\text{ }x_{-}=\max (-x,0).
\end{equation*}
\end{lemma}
\noindent \textit{Proof.} We note $(t,r)$ for $(\chi_0,\chi_1)$
and we define $\hat{g}(\tau ,r)=\mathcal{F}_{t}(a(t,r))(\tau )$,
where $\mathcal{F}_{t}$ is the partial Fourier transform with
respect to $t$. Then, we obtain:
\begin{equation*}
\int\limits_{0}^{\infty }(\int\limits_{\mathbb{R}}e^{i\lambda
tr^{k}}a(t,r)dt)dr=\int\limits_{0}^{\infty }\hat{g}(-\lambda
r^{k},r)dr=\lambda ^{-\frac{1}{k}}\int\limits_{0}^{\infty }\hat{g}%
(-r^{k},\frac{r}{\lambda ^{\frac{1}{k}}})dr.
\end{equation*}
A Taylor expansion with respect to $r$ for $\hat{g}(\tau ,r)$ at
the origin gives:
\begin{equation*}
\hat{g}(-r^{k},\frac{r}{\lambda ^{\frac{1}{k}}})=\sum\limits_{l=0}^{N}\frac{%
\lambda ^{-\frac{l}{k}}}{l!}r^{l}\frac{\partial ^{l}\hat{g}}{\partial r^{l}}%
(-r^{k},0)+\lambda ^{-\frac{N+1}{k}}R_{N+1}(r,\lambda ),
\end{equation*}
where $R_{N+1}(r,\lambda )$ is integrable with respect to $r,$
with $L^{1}$ norm uniformly bounded in $\lambda $. By a new change
of variable we obtain:
\begin{equation*}
\lambda^{-\frac{1}{k}}\int\limits_{0}^{\infty}\hat{g}(-r^{k},\frac{r}{\lambda
^{\frac{1}{k}}})dr = \frac{1}{k} \sum\limits_{l=0}^{N}
\frac{1}{l!}\lambda^{-\frac{1+l}{k}}\int\limits_{-\infty }^{0}
\frac{\partial ^{l}\hat{g}}{\partial r^{l}}
(r,0)|r|^{\frac{l+1-k}{k}}dr+\mathcal{O}(\lambda
^{-\frac{N+1}{k}}).
\end{equation*}
If we introduce $x_{-} =\max (-x,0)$ and
$\mathcal{F}(x_{-}^{\frac{l+1-k}{k}})(r)$ the result follows.
\hfill{$\blacksquare$}
\begin{remark} \rm{A similar result holds for the phase $-\chi
_{0}\chi _{1}^{k}$ if we change $\chi_0$ into $-\chi_0$ and we
then have simply to replace terms $x_{-}$ by
$x_{+}=\mathrm{max}(x,0)$.}
\end{remark}
\section{Proof of the main result.}
We can assume that $z_0$ is a maximum. Lemma \ref{Theo IO 1ere
carte} shows that the first non-zero coefficient is obtained for
$l=2n-1$ (see Remark \ref{shift}) and is given by:
\begin{equation*}
\frac{1}{k}\frac{1}{(2n-1)!}\left\langle
(\mathcal{F}(x_{-}^{\frac{2n-k}{k}})\otimes \delta
_{0}^{(2n-1)}),A(\chi _{0},\chi _{1})\right\rangle
=\frac{1}{k}\int (\mathcal{F}(x_{-}^{\frac{2n-k}{k}})
(\chi_{0})\tilde{A}(\chi _{0},0)d\chi _{0},
\end{equation*}
where, by construction, the new amplitude is:
\begin{equation*}
\tilde{A}(\chi _{0},\chi _{1}) =\int \chi ^{\ast }(a(t,r\theta
)|J\chi ||\mathfrak{p}_{k}(\theta )+rR_{k+1}(\theta
)+tG_{k+1}(t,r\theta )|^{-\frac{2n-1}{k}})d\chi _{2}...d\chi
_{2n}.
\end{equation*}
Now, we use an oscillatory representation of $\tilde{A}(\chi
_{0},0)$ via:
\begin{equation*}
\tilde{A}(\chi _{0},0)=\frac{1}{2\pi}\int\tilde{A}(\chi
_{0},\chi_{1})e^{iz\chi _{1}}dzd\chi_1.
\end{equation*}
We return to the initial coordinates and, since $\chi_{0}=t$, we
obtain:
\begin{equation*}
\tilde{A}(t,0) =\frac{1}{2\pi}\int a(t,r\theta
)|\mathfrak{p}_{k}(\theta )+rR_{k+1}(\theta )+tG_{k+1}(t,r\theta
)|^{-\frac{2n-1}{k}}e^{iz\chi _{1}(t,r,\theta )}dzdrd\theta .
\end{equation*}
Inserting the definition of $\chi _{1}$, we use the change of
variable:
\begin{equation*}
y =z|\mathfrak{p}_{k}(\theta )+rR_{k+1}(\theta
)+tG_{k+1}(t,r\theta )|^{\frac{1}{k}}.
\end{equation*}
Since in $r=0$ we have $G_{k+1}(t,0)=0$, by integration in $(y,r)$
we obtain:
\begin{equation}
\tilde{A}(\chi _{0},0)=\int\limits_{\mathbb{S}^{2n-1}}a(\chi
_{0},0)|\mathfrak{p}_{k}(\theta )|^{-\frac{2n}{k}}d\theta .
\end{equation}
The distribution $\mathcal{F}(x_{-}^{\frac{2n-k}{k}})(\chi _{0})$
of Lemma \ref{Theo IO 1ere carte} acts on this function via:
\begin{equation*}
\left\langle \mathcal{F}(x_{-}^{\frac{2n-k}{k}})(\chi
_{0}),\tilde{A}(\chi _{0},0)\right\rangle =\left\langle
\mathcal{F}(x_{-}^{\frac{2n-k}{k}}) ,a(.,0)\right\rangle
\int\limits_{\mathbb{S}^{2n-1}}|\mathfrak{p}_{k}(\theta
)|^{-\frac{2n}{k}}d\theta .
\end{equation*}
Hence, the top-order contribution to the trace formula is given
by:
\begin{equation*}
\gamma (E_{c},h)=\frac{1}{k}\frac{h^{-n}}{(2\pi
)^{n+1}}h^{\frac{2n}{k}}\left\langle
\mathcal{F}(x_{-}^{\frac{2n-k}{k}}),a(.,0)\right\rangle
\int\limits_{\mathbb{S}^{2n-1}}|\mathfrak{p}_{k}(\theta)|^{-\frac{2n}{k}}d\theta
+\mathcal{O}(h^{\frac{2n+1}{k}-n}).
\end{equation*}
With $a(t,0)=\hat{\varphi}(t)\exp (itp_{1}(z_{0}))$, see
(\ref{demi densite}), by Fourier inversion we have:
\begin{equation*}
\frac{1}{2\pi}\left\langle \mathcal{F}(x_{-}^{\frac{2n-k}{k}}),a(.,0)\right\rangle =\int\limits_{\mathbb{R}%
}\varphi (t+p_{1}(z_{0}))t_{-}^{\frac{2n-k}{k}}dt.
\end{equation*}
Finally, if $z_0$ is a local minimum we must simply replace
$t_{-}$ by $t_{+}=\mathrm{max}(t,0)$ and this complete the proof
of Theorem \ref{Main1}. \hfill{$\blacksquare$}


\begin{thebibliography}{}
\bibitem{BB} R. Balian and C. Bloch, Solution of the Schr\"odinger equation in
terms of classical paths, Annals of Physics \textbf{85} (1974)
514-545.
\bibitem{BPU} R. Brummelhuis, T. Paul and A. Uribe, Spectral estimate near a critical level,
Duke Mathematical Journal \textbf{78} (1995) no. 3, 477-530.
\bibitem{BU} R. Brummelhuis and A. Uribe, A semi-classical trace formula for Schr\"odinger
operators, Communications in Mathematical Physics \textbf{136}
(1991) no. 3, 567-584.
\bibitem{Cam0} B. Camus, Formule des traces semi-classique au
niveau d'une \'energie critique, Th\`ese de l'universit\'e de
Reims (2001).
\bibitem{Cam} B. Camus, A semi-classical trace formula at a non-degenerate critical
level, Journal of functional analysis \textbf{208} (2004), no. 2,
446-481.
\bibitem{C-P} A.M. Charbonnel and G. Popov, A semi-classical trace
formula for several commuting operators, Communications in Partial
Differential Equations \textbf{24} (1999) no. 1-2, 283-323.
\bibitem{CRR} M. Combescure, J. Ralston and D. Robert, A proof of the
Gutzwiller semi-classical trace formula using coherent states
decomposition, Communications in Mathematical Physics \textbf{202}
(1999) no. 2, 463-480.
\bibitem{DUI1} J.J. Duistermaat, Oscillatory integrals Lagrange immersions and unfolding
of singularities, Communications on Pure and Applied Mathematics
\textbf{27} (1974) 207-281.
\bibitem{D-H} J.J. Duistermaat and L. H\"ormander, Fourier Integral
Operators, Acta mathematica \textbf{128} (1972) no. 3-4, 183-269.
\bibitem{GUT} M. Gutzwiller, Periodic orbits and classical quantization conditions,
Journal Math. Phys. \textbf{12} (1971) 343-358.
\bibitem{HOR1} L. H{\"o}rmander, \textit{The analysis of linear partial operators 1,2,3,4} (Springer-Verlag ,1985).
\bibitem{HOR2} L. H{\"o}rmander, \textit{Seminar on singularities of solutions of linear
partial differential equations} (Annals of mathematical studies
91, Princeton University Press, 1979) pages 3-49.
\bibitem{KhD} D. Khuat-Duy, A semi-classical trace formula for Schr\"odinger operators in the
case of a critical level, Th\`ese de l'universit\'e Paris 9
(1996).
\bibitem{KhD1} D. Khuat-Duy, A semi-classical trace formula at a critical level,
Journal of Functional Analysis \textbf{146} (1997) no. 2, 299-351.
\bibitem{PU} T. Paul and A. Uribe, Sur la formule semi-classique des
traces, Comptes Rendus des Sc\'eances de l'Acad\'emie des
Sciences. S\'erie I. \textbf{313} (1991) no. 5, 217-222.
\bibitem{P-P} V. Petkov and G. Popov, Semi-classical trace formula and
clustering of the eigenvalues for Schr\"odinger operators, Annales
de l'Institut Henri Poincar\'e. Physique Th\'eorique. \textbf{68}
(1998) no. 1, 17-83.
\bibitem{[Rob]} D. Robert, \textit{Autour de l'approximation semi-classique}
(Progress in mathematics 68, Birkh\"{a}user Boston, 1987).
\end{thebibliography}
\end{document}